\documentclass[11pt]{article}

\input{epsf}

\usepackage{amsfonts}

\newcommand{\qed}{\hbox{\rule{6pt}{6pt}}}

\newcommand{\Z}{\mathbb{Z}}

\newtheorem{theorem}{Theorem}[section]
\newtheorem{corollary}[theorem]{Corollary}
\newtheorem{lemma}[theorem]{Lemma}
\newtheorem{proposition}[theorem]{Proposition}

\newtheorem{example}[theorem]{Example}

\newtheorem{remark}[theorem]{Remark}

\begin{document}

\title{Extensions of Quandles and Cocycle Knot Invariants }

\author{
J. Scott Carter
\\University of South Alabama \\
Mobile, AL 36688 \\ carter@mathstat.usouthal.edu
\and 
Mohamed Elhamdadi
\\ University of South Florida
\\ Tampa, FL 33620  \\ emohamed@math.usf.edu
\and 
Marina Appiou Nikiforou
\\ University of South Florida
\\ Tampa, FL 33620  \\ mappiou@chuma.cas.usf.edu
\and 
Masahico Saito
\\ University of South Florida
\\ Tampa, FL 33620  \\ saito@math.usf.edu
}

\maketitle

\vspace{10mm}

\begin{abstract}
Quandle cocycles are constructed from extensions of quandles.
The theory is parallel to that of group cohomology and group extensions.
An interpretation of quandle cocycle invariants as obstructions to 
extending knot colorings is given, and is extended to links component-wise.
\end{abstract}

\vspace{10mm}

\section{Introduction}

A quandle is a set with a 
self-distributive binary operation (defined below)
whose study has been partially motivated from knot theory. 
A (co)homology theory was defined in \cite{CJKLS} for quandles,
which is a modification of rack (co)homology defined in \cite{FRS2}. 
State-sum invariants using quandle cocycles as weights are 
defined \cite{CJKLS} and computed for important families
of classical knots and knotted surfaces \cite{CJKS1}.
Quandle homomorphisms and virtual knots are applied to this 
homology theory \cite{betti}.
 The invariants were applied to study 
knots, 
for example, in detecting non-invertible
knotted surfaces \cite{CJKLS}. 
On the other hand, knot diagrams colored by quandles can be used 
to study quandle homology groups. This view point was developed
 in \cite{FRS2,Flower,Greene}  
for rack homology and homotopy, and generalized to quandle homology
in \cite{SSS2}.

In \cite{SanFran}, constructions of extensions of quandles 
using cocycles are given, which are similar to extensions of groups
using group cocycles \cite{Brown}. 
In this paper, we develop methods of constructing cocycles from extensions.
This is the opposite direction of \cite{SanFran}. 
After reviewing the material in Section~2, 
infinite families of extensions of Alexander quandles are given,
and explicit formulas of computing corresponding cocycles are
established in Section~3. We show that these  families of
extensions are non-trivial. 
In Section~4, 
 an interpretation of
the state-sum cocycle invariants as obstructions to 
extending colorings by a quandle to colorings by another quandle
is given. Generalizations to links are also considered.

\section{Quandle and Their Homology Theory}

In this section we review necessary material from the papers mentioned
in the introduction.

A {\it quandle}, $X$, is a set with a binary operation $(a, b) \mapsto a * b$
such that

(I) For any $a \in X$,
$a* a =a$.

(II) For any $a,b \in X$, there is a unique $c \in X$ such that 
$a= c*b$.

(III) 
For any $a,b,c \in X$, we have
$ (a*b)*c=(a*c)*(b*c). $

A {\it rack} is a set with a binary operation that satisfies 
(II) and (III).

Racks and quandles have been studied in, for example, 
\cite{Brieskorn,FR,Joyce,K&P,Matveev}.

The axioms for a quandle correspond respectively to the 
Reidemeister moves of type I, II, and III 
(see 
\cite{FR,K&P}, for example).

A function $f: X \rightarrow  Y$ between quandles
or racks  is a {\it homomorphism}
if $f(a \ast b) = f(a) * f(b)$ 
for any $a, b \in X$.

The following are typical examples of quandles.

\begin{itemize}
\item
A group $X=G$ with
$n$-fold conjugation
as the quandle operation: $a*b=b^{-n} a b^n$.
\item
Any set $X$ with the operation $x*y=x$ for any $x,y \in X$ is
a quandle called the {\it trivial} quandle.
The trivial quandle of $n$ elements is denoted by $T_n$.
\item
Let $n$ be a positive integer.
For elements  $i, j \in \{ 0, 1, \ldots , n-1 \}$, define
$i\ast j \equiv 2j-i \pmod{n}$.
Then $\ast$ defines a quandle
structure  called the {\it dihedral quandle},
  $R_n$.
This set can be identified with  the
set of reflections of a regular $n$-gon
  with conjugation
as the quandle operation.
\item
Any $\Lambda (={\Z }[T, T^{-1}])$-module $M$
is a quandle with
$a*b=Ta+(1-T)b$, $a,b \in M$, called an {\it  Alexander  quandle}.
Furthermore for a positive integer
$n$, a {\it mod-$n$ Alexander  quandle}
${\Z }_n[T, T^{-1}]/(h(T))$
is a quandle
for
a Laurent polynomial $h(T)$.
The mod-$n$ Alexander quandle is finite
if the coefficients of the
highest and lowest degree terms
of $h$
  are units in $\Z_n$. 
\end{itemize}

\begin{figure}
\begin{center}
\mbox{
\epsfxsize=3in
\epsfbox{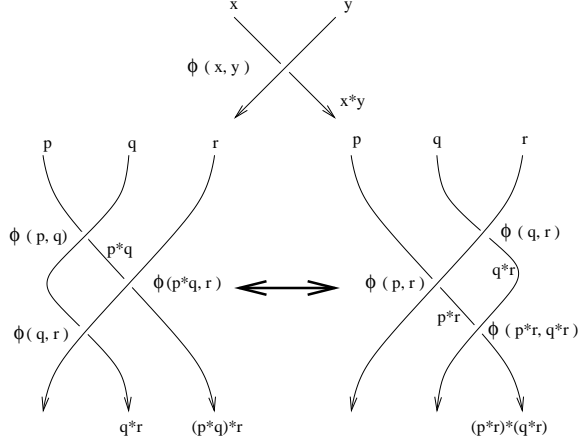} 
}
\end{center}
\caption{ Type III move and the quandle identity  }
\label{2cocy} 
\end{figure}

\bigskip

 Let $C_n^{\rm R}(X)$ be the free 
abelian group generated by
$n$-tuples $(x_1, \dots, x_n)$ of elements of a quandle $X$. Define a
homomorphism
$\partial_{n}: C_{n}^{\rm R}(X) \to C_{n-1}^{\rm R}(X)$ by \begin{eqnarray}
\lefteqn{
\partial_{n}(x_1, x_2, \dots, x_n) } \nonumber \\ && =
\sum_{i=2}^{n} (-1)^{i}\left[ (x_1, x_2, \dots, x_{i-1}, x_{i+1},\dots, x_n) \right.
\nonumber \\
&&
- \left. (x_1 \ast x_i, x_2 \ast x_i, \dots, x_{i-1}\ast x_i, x_{i+1}, \dots, x_n) \right]
\end{eqnarray}
for $n \geq 2$ 
and $\partial_n=0$ for 
$n \leq 1$. 
 Then
$C_\ast^{\rm R}(X)
= \{C_n^{\rm R}(X), \partial_n \}$ is a chain complex.

Let $C_n^{\rm D}(X)$ be the subset of $C_n^{\rm R}(X)$ generated
by $n$-tuples $(x_1, \dots, x_n)$
with $x_{i}=x_{i+1}$ for some $i \in \{1, \dots,n-1\}$ if $n \geq 2$;
otherwise let $C_n^{\rm D}(X)=0$. If $X$ is a quandle, then
$\partial_n(C_n^{\rm D}(X)) \subset C_{n-1}^{\rm D}(X)$ and
$C_\ast^{\rm D}(X) = \{ C_n^{\rm D}(X), \partial_n \}$ is a sub-complex of
$C_\ast^{\rm
R}(X)$. Put $C_n^{\rm Q}(X) = C_n^{\rm R}(X)/ C_n^{\rm D}(X)$ and 
$C_\ast^{\rm Q}(X) = \{ C_n^{\rm Q}(X), \partial'_n \}$,
where $\partial'_n$ is the induced homomorphism.
Henceforth, all boundary maps will be denoted by $\partial_n$.

For an abelian group $G$, define the chain and cochain complexes
\begin{eqnarray}
C_\ast^{\rm W}(X;G) = C_\ast^{\rm W}(X) \otimes G, \quad && \partial =
\partial \otimes {\rm id}; \\ C^\ast_{\rm W}(X;G) = {\rm Hom}(C_\ast^{\rm
W}(X), G), \quad
&& \delta= {\rm Hom}(\partial, {\rm id})
\end{eqnarray}
in the usual way, where ${\rm W}$ 
 $={\rm D}$, ${\rm R}$, ${\rm Q}$.

The $n$\/th {\it quandle homology group\/}  and the $n$\/th
{\it quandle cohomology group\/ } \cite{CJKLS} of a quandle $X$ with coefficient group $G$ are
\begin{eqnarray}
H_n^{\rm Q}(X; G) 
 = H_{n}(C_\ast^{\rm Q}(X;G)), \quad
H^n_{\rm Q}(X; G) 
 = H^{n}(C^\ast_{\rm Q}(X;G)). \end{eqnarray}

\bigskip

Let a classical knot diagram be given. 
The co-orientation is a family of normal vectors to the knot diagram 
such that the pair (orientation, co-orientation) matches
the given (right-handed, or counterclockwise) orientation of the plane.
At a crossing, 
if the pair of the co-orientation 
 of the 
over-arc and  that of the under-arc
matches the (right-hand) orientation of the plane, then the 
crossing is called {\it positive}; otherwise it is {\it negative}. 
Crossings in Fig.~\ref{2cocy} are positive by convention.

A  {\it coloring}   
of an oriented  classical knot diagram is a
function ${\mathcal C} : R \rightarrow X$, where $X$ is a fixed 
quandle
and $R$ is the set of over-arcs in the diagram,
satisfying the  condition
depicted 
in the top
of Fig.~\ref{2cocy}. 
In the figure, a 
crossing with
over-arc, $r$, has color ${\mathcal C}(r)= y \in X$. 
The under-arcs are called $r_1$ and $r_2$ from top to bottom;
the normal (co-orientation) of the over-arc $r$ points from $r_1$ to $r_2$.
Then it is required that 
${\mathcal C}(r_1)= x$ and ${\mathcal C}(r_2)=x*y$.

Note that locally the colors do not depend on the 
orientation of the under-arc.
The quandle element ${\mathcal C}(r)$ assigned to an arc $r$ by a coloring 
 ${\mathcal C}$ is called a {\it color} of the arc. 
This definition of colorings on knot diagrams has been known, see 
\cite{FR,FoxTrip} for example. 
Henceforth, all the quandles that are used to color diagrams will be finite.

In Fig.~\ref{2cocy} bottom, the relation between Redemeister type III move
and a quandle axiom (self-distributivity) is indicated. 
In particular, the colors of the bottom right segments before and after
the move correspond to the self-distributivity. 

\bigskip

A {\it (Boltzmann) weight}, $B(\tau, {\cal C})$,
at a  crossing $\tau$ is defined as follows.
Let  ${\cal  C}$ 
denote a coloring.
Let $r$ be the over-arc at $\tau$, and $r_1$, $r_2$ be 
under-arcs such that the 
normal
to $r$ points from $r_1$ to $r_2$.
Let $x={\cal C}(r_1)$ and $y={\cal C}(r)$. 
 Pick a 
quandle 
2-cocycle 
$\phi \in  Z^2(X; A)$.
Then define $B(\tau, {\cal C})= \phi(x,y)^{\epsilon (\tau)}$,
where 
$\epsilon (\tau)= 1$ or $-1$, if  the sign of $\tau$ 
is positive or negative, respectively.

The {\it partition function}, or a {\it state-sum}, 
is the expression 
$$
\sum_{{\cal C}}  \prod_{\tau}  B( \tau, {\cal C}).
$$
The product is taken over all crossings of the given diagram,
and the sum is taken over all possible colorings.
The values of the partition function 
are  taken to be in  the group ring ${\Z }[A]$ where $A$ is the coefficient 
group 
written multiplicatively. 
This is proved \cite{CJKLS} to be a knot invariant, called
the {\it (quandle) cocycle invariant}. 
\section{Extensions of quandles by $2$-cocycles}

In \cite{SanFran}, for a quandle $X$, an abelian group $A$, and a $2$-cocycle 
 $\phi \in Z^2_{\rm Q}(X; A)$, the {\it abelian extension} 
$E=E(X,A,\phi)$ was defined as $A \times X$ as a set, 
with the quandle operation 
 $(a_1, x_1) * (a_2, x_2) = (a_1 \phi(x_1, x_2), x_1 * x_2)$.
The following lemma is the converse of the fact proved in \cite{SanFran} that 
$E(X,A,\phi)$ is a quandle.

\begin{lemma} \label{cocylemma}
Let $X$, $E$ be finite quandles, and $A$ be a finite abelian group
written multiplicatively. 
Suppose there exists a bijection 
$f: E \rightarrow A \times X$ with the following property.
There 
exists 
a function $\phi: X \times X \rightarrow A$ such that
for any $e_i \in E$ ($i=1,2$), 
if $f(e_i)=(a_i, x_i)$, then 
$f(e_1 *  e_2) = (a_1 \phi(x_1, x_2) , x_1 * x_2 )$. 
Then $\phi \in Z^2_{\rm Q}(X; A)$. 
\end{lemma}  
 {\it Proof.\/} 
For any $x \in X$ and $a \in A$,
 there is $e \in E$ such that   $f(e)=(a,x)$, and 
$$ (a,x)=f(e)=f(e*e)=(a \phi(x, x), x) , $$
so that we have $\phi(x, x)=1$ for any $x \in X$.

By identifying $A \times X$ with $E$ by $f$, the quandle operation
$*$ on  $A \times X$ is defined, for any  $(a_i, x_i)$ ($i=1,2$), by
 $$(a_1, x_1) * (a_2, x_2) = (a_1 \phi(x_1, x_2), x_1 * x_2).$$
Since  $A \times X$ is a quandle isomorphic to $E$ under this $*$, we have 
\begin{eqnarray*}
\lefteqn{
  [ (a_1, x_1) * (a_2, x_2)] * (a_3, x_3) } \\
 &=&
(a_1 \phi(x_1, x_2), x_1 * x_2) * (a_3, x_3) \\
&=& (a_1 \phi(x_1, x_2) \phi(x_1 * x_2,  x_3) ,  (x_1 * x_2)* x_3 ) , 
\end{eqnarray*} 
and
\begin{eqnarray*}
\lefteqn{
 [ (a_1, x_1) * (a_3, x_3)] *[ (a_2, x_2) * (a_3, x_3)]} \\
 &=& 
 (a_1 \phi(x_1, x_3), x_1 * x_3) *  (a_2 \phi(x_2, x_3), x_2 * x_3) \\
&=& (a_1 \phi(x_1, x_3) \phi( x_1 * x_3,   x_2 * x_3), 
(x_1 * x_3) *(x_2 * x_3) )
\end{eqnarray*}
are equal for any  $(a_i, x_i)$ ($i=1,2,3$).
Hence $\phi$ satisfies the $2$-cocycle condition.
$\qed$

This lemma implies that under the same assumption  
we have $E=E(X,A,\phi)$,  where  $\phi \in Z^2_{\rm Q}(X; A)$.
Next we identify such examples.

\begin{theorem} \label{firstextthm}
For any  positive integers $q$ and $m$, 
$E=\Z _{q^{m+1}} [T, T^{-1}] /  (T -1 +q) $ is an abelian extension 
$E=E( \Z _{q^{m}} [T, T^{-1}] /  (T -1 +q) ,  \Z_q,  \phi)$
of $X= \Z _{q^{m}} [T, T^{-1}] /  (T -1 +q)$ for some cocycle
 $\phi \in  Z^2_{\rm Q}( X;  \Z_q)$. 
\end{theorem}
{\it Proof.\/}
Represent elements of $\Z _{q^{m+1}}$ by $\{ 0, 1, \ldots, q^{m+1} -1 \}$
and express them 
in their
 $q^{m+1}$-ary expansion:
$$A = A_{m} q^m + \ldots + A_1 q + A_0 \in \Z _{q^{m+1}}, $$
where $0 \leq A_j < q$, $j=0, \ldots, m$. 
With this convention, $A_j$'s are uniquely determined integers. 
Define  $f: E \rightarrow \Z_q \times X$ by
$$f(A)=(A_m \  (\mbox{mod}\ q), \  \overline{A} \  (\mbox{mod}\ q^{m-1})
), $$
where 
$ \overline{A}=  \sum_{j=0}^{m-1} A_{j} q^{j} . $
Then for $A, B  \in \Z _{q^{m+1}}$, the quandle operation is computed 
in $\Z_{q^{m+1} }$ by
\begin{eqnarray*} 
 A*B & = & TA + (1-T)B \\
 & = & (1-q) ( A_{m} q^m + \ldots + A_1 q + A_0 ) \\
 & & + q( B_{m} q^m + \ldots + B_1 q + B_0 ) \\
 & = & ( A_{m} - A_{m-1} + B_{m-1}  ) q^m
 + ( A_{m-1} - A_{m-2}+B_{m-2}  ) q^{m-1} \\
 & & + \ldots + ( A_1 - A_0+B_0 ) q + A_0 \\
 & = &  ( A_{m} - A_{m-1} + B_{m-1} ) q^m +
 \sum_{j=0}^{m-1} (A_{j} - A_{j-1}+B_{j-1}) q^j ,
\end{eqnarray*} 
where $A_{-1}, B_{-1} $ are understood to be zeros in the last summation.
Define a set-theoretic section
 $s: {\Z }_{q^m} \rightarrow {\Z }_{q^{m+1}}$ by 
$$s \left( \sum_{j=0}^{m-1} X_j q^j \right) =  0 \cdot q^{m}+ 
\sum_{j=0}^{m-1} X_j q^j.$$
For $X, Y\in {\Z }_{q^m}$ define 
$$\phi(X,Y) =  [s(X) * s(Y) - s(X*Y) ] / q^m  
 \in {\Z }_q$$
 where division by  $q^m$ 
means to consider these elements as integers,
divide by  $q^m$ 
and compute the residue class modulo $q$. 
Note that $\overline{s(X)} *\overline{s(Y) }  = \overline{s(X*Y)}.$
Hence,
$s(X) * s(Y) - s(X*Y)$ is divisible by  $q^m$. 
Then we have 
 $$f(A*B) = ( A_m + \phi(\overline{A}, \overline{B}),
\overline{A}*\overline{B}). $$
So that $f$ yields an isomorphism 
$$\Z_{q^{m+1}} [T, T^{-1}] /  (T -1 +q) \rightarrow
 E( \Z _{q^{m}} [T, T^{-1}] /  (T -1 +q) ,  \Z_q,  \phi). \quad  
\qed$$

\begin{theorem} \label{secondextthm}
For any positive integer $q$ and $m$, the quandle 
$E=\Z_q [T, T^{-1} ] / (1-T)^{m+1} $ is an abelian extension
of $X=\Z_q [T, T^{-1} ] / (1-T)^{m} $ over $\Z_q$: 
$E=E(X, \Z_q, \phi)$, for some $\phi \in Z^2_{\rm Q}(X; \Z_q)$. 
\end{theorem}
{\it Proof.\/} 
Represent elements of $E$ by 
$A=A_m (1-T)^m + \ldots + A_1(1- T) + A_0$, where
$A_j \in \Z_q $, $j=0, \ldots, m$.
Define  $f: E \rightarrow \Z_q \times X$ by
$$f(A)=(A_m \  (\mbox{mod}\ (1-T) ), 
 \  \overline{A} \  (\mbox{mod}\ (1-T)^{m-1}) ), $$
where 
$ \overline{A}=  \sum_{j=0}^{m-1} A_{j} (1-T)^{j} . $
Then for $A, B  \in E$, the quandle operation is computed 
 by
\begin{eqnarray*} 
 A*B & = & TA + (1-T)B \\
 & = & [1-(1-T)]  ( A_{m} (1-T)^m + \ldots + A_1 (1-T) + A_0 ) \\
 & & + (1-T) ( B_{m} (1-T)^m + \ldots + B_1 (1-T) + B_0 ) \\
 & = & ( A_{m}- A_{m-1}+ B_{m-1}  ) (1-T)^m
 + ( A_{m-1} - A_{m-2}+B_{m-2}  ) (1-T)^{m-1} \\
 & & + \ldots + ( A_1 - A_0+B_0 ) (1-T) + A_0 \\
 & = &  ( A_{m}- A_{m-1}+  B_{m-1}  ) (1-T)^m +
 \sum_{j=0}^{m-1} (A_{j} - A_{j-1}+B_{j-1}) (1-T)^j , 
\end{eqnarray*} 
where $A_{-1}, B_{-1} $ are understood to be zeros in the last summation,
and the coefficients are in $\Z_q$.
 Note that 
in $\Z_q[T, T^{-1}]/{(1-T)^{m}}$, we have 
$$\overline{A}*\overline{B} = 
\sum_{j=0}^{m-1} (A_{j} - A_{j-1}+B_{j-1}) (1-T)^j . $$
Hence we have
$$f ( A*B )= (   A_{m}- A_{m-1} + B_{m-1} , 
 \sum_{j=0}^{m-1} (A_{j} - A_{j-1}+B_{j-1}) (1-T)^j ) \in \Z_q \times X. $$
Then we have
$$ f(A*B)=( A_m + 
\phi(\overline{A}, \overline{B}), \overline{A}*\overline{B}),$$
where $\phi(\overline{A}, \overline{B})= B_{m-1}- A_{m-1} $.
Hence $f$ yields an isomorphism 
$$\Z_{q} [ T, T^{-1}] / (1-T)^{m+1}\rightarrow
 E( \Z _{q} [T, T^{-1} ] / (1-T)^{m} ,  \Z_q,  \phi).$$

The cocycle $\phi$ has a similar description to the one in
 Theorem~\ref{firstextthm}. Let 
$$s: \Z_{q} [ T, T^{-1}] / (1-T)^{m} \rightarrow \Z_{q} [ T, T^{-1}] / (1-T)^{m+1}$$ be a set-theoretic section defined by 
$$s\left( \sum_{j=0}^{m-1} A_j (1-T)^j \quad \mbox{mod}\ (1-T)^{m} \right)
=  \sum_{j=0}^{m-1} A_j (1-T)^j \quad  \mbox{mod}\ (1-T)^{m+1}. $$
Then we have 
$\overline{ s(X)} * \overline{ s(Y)} =\overline{ s(X *Y)}$ for any
$X, Y \in   \Z_{q} [ T, T^{-1}] / (1-T)^{m}$.  
So that 
$[  s(X) *  s(Y)- s(X *Y)]$ is divisible by $(1-T)^m$, 
and we have 
$$ \phi(\overline{A}, \overline{B})
=  [  s(X) *  s(Y)- s(X *Y)]/ (1-T)^m \in \Z_q. \quad \qed$$

\begin{example} {\rm 

\begin{enumerate}

\item
Consider the case $q=2$, $m=2$ in Theorem~\ref{firstextthm}.
In this case 
$$\Z_4[T, T^{-1}]/ (T+1)=R_4, \quad \mbox{and} $$ 
 $$\Z_8[T, T^{-1}]/ (T+1)=R_8=E(R_4, \Z_2, \phi)$$ for some 
$\phi \in Z^2_{\rm Q}(R_4; \Z_2)$.
We  obtain an explicit formula for this cocycle
$\phi$ by computation: 
$$\phi = \chi_{0,2} + \chi_{0,3} + \chi_{1,0} + \chi_{1,3}
+ \chi_{2,0} + \chi_{2,3} + \chi_{3,0} + \chi_{3,1}   , $$
where 
$$\chi_{a,b} (x,y) = \left\{ \begin{array}{ll} 1 & {\mbox{\rm if }} \
(x,y)=(a,b), \\
0 & {\mbox{\rm if }} \
(x,y)\not=(a,b) \end{array}\right.$$
denotes the characteristic function.

\item In case $m=1$ and $q=3$, the cocycle constructed is of the form
$$ \phi= \chi_{0,1} + \chi_{1,2} + \chi_{2,0} +2(\chi_{0,2} + \chi_{1,0} + \chi_{2,1}).$$

\item In case $m=2$ and $q=3$, the cocycle is

$$
\begin{array}{lll}
 \phi & =  &\chi_{0,3} +\chi_{0,4}+ \chi_{0,5}+
2 \chi_{0,6}+2\chi_{0,7}+2\chi_{0,8} \\
&+& 2 \chi_{1,0} + \chi_{1,4}+\chi_{1,5}+\chi_{1,6}+2\chi_{1,7}+2\chi_{1,8}\\
&+& 2 \chi_{2,0}+2 \chi_{2,1}+ \chi_{2,5}+ \chi_{2,6}+ \chi_{2,7}
+ 2 \chi_{2,8} \\
&+& 2 \chi_{3,0} + 2 \chi_{3,1} +  \chi_{3,5} +  \chi_{3,6} + 
 \chi_{3,7} + 2 \chi_{3,8} \\
&+&  2\chi_{4,0} + 2 \chi_{4,1}+2 \chi_{4,2}+   \chi_{4,6}
+  \chi_{4,7}+ \chi_{4,8}\\
&+&  \chi_{5,0}+2 \chi_{5,1}+2 \chi_{5,2}+2 \chi_{5,3}
+ \chi_{5,7}+ \chi_{5,8} \\
&+&  \chi_{6,0}+2 \chi_{6,1}+2 \chi_{6,2}+2 \chi_{6,3}+
 \chi_{6,7}+ \chi_{6,8} \\
&+&  \chi_{7,0}+ \chi_{7,1}+2 \chi_{7,2}+2 \chi_{7,3}+2 \chi_{7,4}
+  \chi_{7,8} \\
&+&  \chi_{8,0}+ \chi_{8,1}+ \chi_{8,2}+2 \chi_{8,3}+2 \chi_{8,4}
+2 \chi_{8,5} 
\end{array}
$$

\item
Consider the case $q=2$ and $m=2$ in Theorem~\ref{secondextthm}. 
The quandle $\Z_2[T,T^{-1}]/(1-T)^2$ is isomorphic to $R_4$ by 
 the correspondence  
$0 \leftrightarrow 0 (1-T) + 0 $, 
 $1 \leftrightarrow 0 (1-T) + 1 $, 
 $2 \leftrightarrow 1 (1-T) + 0 $, and 
 $3 \leftrightarrow 1 (1-T) + 1 $.
This is a special case of the isomorphism
$$\Z_{n} [T, T^{-1} ] / (1-T )^2  \cong \Z_{n^2} [T, T^{-1} ] /(T-(kn+1)) 
  \quad  \mbox{if} \quad (n,k)=1  $$
given in \cite{LithNel}. 
Then the quandle $\Z_2[T,T^{-1}]/(1-T)^3$ is an abelian extension
$E( R_4 ; \Z_2, \phi ')$ for some 
$\phi'  \in Z^2_{\rm Q} (R_4;\Z_2)$. 
 Then the cocycle $\phi ' (\overline{A}, \overline{B})=B_1 - A_1$ is 
$1$ if and only if the pair  $(\overline{A}, \overline{B})$
has distinct coefficients for $(1-T)$, and we obtain 
$$ \phi'= \chi_{0,2} + \chi_{2,0} + \chi_{1,2} + \chi_{2,1}
 + \chi_{0,3} + \chi_{3,0} + \chi_{1,3} + \chi_{3,1} . $$
The cocycles $\phi_0=\chi_{2,1}+\chi_{2,3}$, 
$\phi_1=\chi_{1,0}+\chi_{1,2}$, and the above $\phi$ are 
linearly independent 
(evaluate on the cycles defined in 
Remark~\ref{rank3cor} below), 
and $\phi' = \phi + \phi_0 + \phi_1$. 
\end{enumerate}
} \end{example}

We recall from \cite{SanFran} that 
two surjective homomorphisms of quandles $\pi_j: E_j \rightarrow X$,
$j=1,2$,
are called {\it equivalent } if there is
a quandle isomorphism $f: E_1 \rightarrow E_2$
such that $\pi_1=\pi_2 f$.
In particular, abelian extensions define surjective homomorphisms
$E(X,A,\phi)=A \times X \rightarrow X$ defined by the projection
onto the second factor. 
It was proved in \cite{SanFran} that two abelian extensions
$E(X,A,\phi)$ and $E(X,A,\phi')$ 
are equivalent in the above sense if and only if 
$\phi$
is cohomologous to $\phi'$.

\begin{proposition}
The abelian extensions 
$$\Z_{q^{m+1}}[T, T^{-1}] / (T-1+q)= E(\Z_{q^m}[T, T^{-1}] / (T-1+q), \Z_q, \phi )  , $$ 
$$ \Z_{q}[T, T^{-1}] / (1-T)^{m+1} = E(\Z_{q}[T, T^{-1}] / (1-T)^m,  \Z_q, 
\phi' ) $$ 
are not trivial, i.e., not  product quandles.
\end{proposition}
{\it Proof.\/} 
Direct computations show that 
the chains
$$ c =(0,1) + (q, q^{m-1}+q-1) \quad \in Z_2^{\rm Q}(X;\Z_q) \quad \mbox{and} $$
$$c'=(0,1) + ( 1-T, (1-T)^{m-1}  + (1-T) -1 ) \; \in Z_2^{\rm Q}(X;\Z_q) $$
are  cycles for $X=\Z_{q^m}[T, T^{-1}] / (T-1+q)$
and $X=\Z_{q}[T, T^{-1}]/(1-T)^m$, respectively.

Then it is computed that $\phi(c)=1$ and $\phi'(c')=1$,
and hence $\phi$ and $\phi'$ are not coboundaries, and the result follows.
$\qed$

\begin{remark} \label{prodremark} {\rm 
We remark here on consequences on dihedral quandles we derive from
the above results. 
The quandle structure of a dihedral quandle $R_n$
is defined using ring structure 
of $\Z_n$. The product quandle $R_m \times R_n$ is defined 
by component-wise operation, so that it is defined from the ring 
structure of $\Z_m \times \Z_n$ as well. 
Consequently, two quandles $R_m \times R_n$ and $R_{mn}$ are isomorphic 
if  $\Z_m \times \Z_n$ and  $\Z_{mn}$ are isomorphic 
as rings.
Hence if $n=p_1^{e_1} \ldots p_k^{e_k}$ is the prime decomposition, 
then $R_n$ is isomorphic to $R_{p_1^{e_1}} \times \ldots \times R_{p_k^{e_k}}$.
For $p=2$, the result of
this section shows that $R_{p^e}$ is described succesively as an extension
of $R_{p^{e-1}}$. 
}  \end{remark}

The following lemma follows from definitions.

\begin{lemma} \label{prodlemma}
Let $X, Y$ be quandles and  $A$ be an abelian group. 
If $E$ is an abelian extension of $X$ for  $\phi \in Z^2_{\rm Q}(X;A)$:
$E=E(X,A, \phi)$, then $E \times Y$ is an abelian extension
of $X \times Y$ for $p^{\#} \phi \in  Z^2_{\rm Q}(X \times Y ;A)$:
$E \times Y= E(X \times Y, A, p^{\#} \phi ) $, 
where $p: X \times Y \rightarrow X$ is the projection to the first factor.
\end{lemma}

\begin{corollary} \label{dihedlemma}
For any  positive integer $n$,  
 $E=R_{4n}$ is an abelian extension $E=E(R_{2n}, \Z_2, \phi)$
of $X=R_{2n}$ for some cocycle $\phi \in  Z^2_{\rm Q}(R_{2n}; \Z_2)$.
\end{corollary}
{\it Proof.\/} 
Let $2n=2^m k$ for an odd integer $k$. 
Then $R_{2n} \cong R_{2^m} \times R_k$ by Remark~\ref{prodremark},
and by Lemma~\ref{prodlemma},  $R_{4n} \cong R_{2^{m+1} } \times R_k$
is an abelian extension of $R_{2n}$ if  $R_{2^{m+1} }$
is an abelian extension of $R_{2^m}$. This follows from 
Theorem~\ref{firstextthm} since $R_{2^m} \cong \Z_{2^m} [T, T^{-1} ]/(T+1)$.
$\qed$

\bigskip

\begin{remark} 
\label{rank3cor} {\rm 
By Lemma~\ref{cocylemma} and Cor.~\ref{dihedlemma}, 
there is a cocycle $\phi \in Z^2_{\rm Q}(R_{4n};\Z_2)$ such that 
$R_{8n}$ is isomorphic to $E(R_{4n}, \Z_2, \phi)$.

Let $\phi_{0, 1}, \phi_{1,0} \in Z^2_{\rm Q}(R_{4n};\Z_2)$
be cocycles defined by 
$$\phi_{0,1}=p^{\#}( \chi_{0,1} + \chi_{0,3} ) , \quad 
\mbox{and} \quad 
\phi_{1,0}=p^{\#}( \chi_{1,0} + \chi_{1,2} ), $$ 
respectively,
where $p: R_{4n} \rightarrow R_4$ is a natural map
$p(x \; \mbox{mod}\ (4n) ) = x \; \mbox{mod}\ (4)$. 
Here, it is known \cite{CJKLS} that  
$$ \chi_{0,1} + \chi_{0,3} , \quad 
\mbox{and} \quad 
\chi_{1,0} + \chi_{1,2}$$
are cocycles in $Z^2_{\rm Q}(R_4; \Z_2)$. 
It is directly computed that 
$$c_{0,1}=(0,1)+(2,1), \quad c_{1,0}=(1,0)+(4n-1,0), \quad  
c_{0,1}'=(0,1)+(2, 2n+1) \in Z_2^{\rm Q}(R_{4n};\Z_2)$$ 
are  cycles.
Then we have 
$$ \begin{array}{lllllllll} 
\phi_{0,1} (c_{0,1})&= & 1, & \phi_{0,1} (c_{1,0}) &= & 0, & \phi_{0,1} (c_{0,1}') &= & 1, \\
\phi_{1,0} (c_{0,1}) &=  &0, & \phi_{1,0} (c_{1,0}) &=  &1, & \phi_{1,0} (c_{0,1}') &=  &0, \\
\phi \ \    (c_{0,1}) &=  &0, & \phi\ \    (c_{1,0}) &= &0, & \phi \ \   (c_{0,1}') &= &1,
\end{array}
$$
Hence we see that the cocycles $\phi_{0,1}, \phi_{1,0}$, and $\phi$
are linearly independent. 

In \cite{Mochizuki,LithNel}, ranks of homology groups are determined for
certain families of  quandles. 
} \end{remark}

\section{Cocycle knot invariants as obstructions to extending colorings}

Let $K$ be a knot, and denote by $\Phi_K(X, \phi)$  the state-sum
invariant of $K$ with respect to a quandle $X$ and a cocycle 
$\phi \in Z^2_{\rm Q}(X; A)$, where $A$ is an abelian group. 
Let $E= E(X,A,\phi)$ be the abelian extension of $X$ by $\phi$.
We characterize when  the state-sum invariant defined from cocycles 
is non-trivial, if the 
cocycles used are 
 those defined from 
abelian extensions. For characterizations on the triviality of colorings,
see \cite{Inoue}.

Let ${\cal C}$ be a coloring of $K$ by $X$. Let $E$ be
a quandle with a surjective homomorphism 
$p: E \rightarrow X$. In this case $E$ is called an {\it extension} of $X$.
If there is a coloring ${\cal C}'$ of $K$ by $E$ such that for 
every arc $a$ of $K$, it holds that $p ({\cal C}'(a))={\cal C}(a)$, 
then  ${\cal C}'$ is called an  {\it extension} of  ${\cal C}$.

\begin{theorem} \label{knotcolorthm}
Let $C_0(K, X)$ be the constant term  
(a positive integer)
of $\Phi_K(X, \phi) $, and  $C(K, X)$ be the number of all  colorings
of $K$ by $X$. Then the number of colorings of $K$ by $X$
that extend to  colorings of $K$ by $E(X,A,\phi)$
is equal to $C_0(K, X)$, and the number of colorings  that do not 
extend is  $C(K, X)- C_0(K, X)$. 
\end{theorem}
{\it Proof.\/}
 Let ${\cal C}$ be a coloring whose contribution to 
$\Phi_K(X, \phi) $ is  $1$. 
Fix this coloring in what follows.
Pick a base point $b_0$ on a knot diagram of $K$.
Let $x \in X$ be the color on the arc $\alpha_0$ containing $b_0$.
Let $\alpha_i$, $i=1, \ldots, n$, be the set of arcs that appear 
in this order when the diagram $K$ is traced in the given orientation of $K$,
starting from $b_0$.
Pick an element $a \in A$ and give a color 
$(a, x)$ on  $\alpha_0$,
so that we define a coloring ${\cal C}'$ by $E$ on $\alpha_0$
by ${\cal C}'(\alpha_0)=(a, x)\in E$. 
We try to extend it to the entire
diagram by traveling the diagram from $b_0$ along the arcs $\alpha_i$,
 $i=1, \ldots, n$, in this order, by induction. 

Suppose  ${\cal C}'(\alpha_i)$ is defined for $0 \leq i < k$. 
Define  ${\cal C}'(\alpha_{k+1})$ as follows.
Suppose that the crossing $\tau_k$ separating $\alpha_k$ and $\alpha_{k+1}$
is positive, and the over-arc at  $\tau_k$ is $\alpha_j$. 
Let  ${\cal C}'(\alpha_{k})=(a, x)$ and ${\cal C}(\alpha_j)=y \in X$.
Then we have ${\cal C}(\alpha_{k+1})=x*y \in X$. 
Define  ${\cal C}'(\alpha_{k+1})=(a  \phi(x, y), x*y)$ in this case.

 Suppose that the crossing $\tau_k$ is negative.
Let  ${\cal C}'(\alpha_{k})=(a, x)$ and ${\cal C}(\alpha_j)=y \in X$.
Then if  ${\cal C}(\alpha_{k+1})=z$, then we have $z*y=x$.
Define  ${\cal C}'(\alpha_{k+1})=(a  \phi(z, y)^{-1}, z)$ in this case.

Define ${\cal C}'(\alpha_{i})$ inductively for all $i=0, \ldots, n$.
Regard $\alpha_0$ as $\alpha_{n+1}$, and repeat the 
above 
construction at the last crossing $\tau_n$ to come back to $\alpha_0$.
By the construction we have 
 ${\cal C}'(\alpha_{n+1})=(a  \prod_{\tau} B(\tau, {\cal C}),
{\cal C}(\alpha_0) )$, where  $\prod_{\tau} B(\tau, {\cal C})$
is the state-sum contribution (the product of Boltzmann weights over all crossings) of  ${\cal C}$.
This   contribution is  equal to $1$ by the assumption that 
$\prod_{\tau} B(\tau, {\cal C})=1$,
and we have a well-defined coloring ${\cal C}'$. 
Hence this color extends to $E(X,A,\phi)$.

Conversely, if a coloring
${\cal C}$ by $X$ extends to a coloring by $E(X,A,\phi)$, 
then from the above argument, we have that  
$(a, x)=(a \prod_{\tau} B(\tau, {\cal C}), x )$,
if $(a, x)$ is the color on the base point $b_0$.
Hence $ \prod_{\tau} B(\tau, {\cal C})=1$.
$\qed$

Thus the non-trivial value of  $\Phi_K(X, \phi) $ is the obstruction
to extending colorings of $K$ by $X$ to $E(X,A,\phi)$, in the following sense: 
there is a coloring ${\cal C}$ of $K$ by $X$ which does not extend to 
a coloring by $E(X,A,\phi)$,
 if and only if $\Phi_K(X, \phi) $ is not a positive integer.

 \begin{corollary} \label{z2case}
For $A=\Z_2$,  $\Phi_K(X, \phi) = a + b t$ (where $t$ is the variable,
the  generator of $\Z_2$) is determined 
by the number of colorings with $X$ and $E$: $a$ is the number of colorings of $X$ that extend to colorings by $E$, and $b$ is the number of those that do not.
\end{corollary}

\begin{example} {\rm 
For $X=\Z_2[T, T^{-1}] / (T^2 + T + 1 )$ 
we have a cocycle  $\phi=\prod_{a, b \neq T} \chi_{(a, b)} \in Z^2_{\rm Q}(X; \Z_2)$ (see \cite{CJKLS}),
let $E=E(X, \Z_2, \phi)$. Then $\Phi_K(S_4, \phi) = a + b t$
has the above characterization.

Specifically, in \cite{CJKS1}, it was computed that among knots
 in the table up to 9 crossings, the state-sum invariant with the above
quandle and the cocycle takes the value $4 + 12 t$ for the knots
$$ 3_1, 4_1, 7_2, 7_3, 8_1, 8_4, 8_{11}, 8_{13}, 9_1, 9_{12}, 9_{13}, 9_{14}, 9_{21}, 9_{23}, 9_{35}, 9_{37}, $$
and the value $16 + 48t$ for $8_{18}, 9_{40}$. 
Hence for the knots in the former list, the number of colorings by $X$
which  extend to those by $E$ is $4$ (trivial colorings, by a single color),
and those that do not extend is $12$ (all non-trivial colorings do not extend).
For the knots in the latter list, there are $16$ colorings that extend,
and $48$ colorings that do not. 
 }\end{example}

We generalize the invariants to links component-wise.
Let $L=K_1 \cup \cdots \cup K_n$ be a link diagram. 
The following generalization of the state-sum invariant,
which follows from Reidemeister moves,  is suggested 
 by 
Rourke and Sanderson in 
personal communications.

\begin{lemma} 
Let $UK_i$, $i=1, \cdots, n$, denote the set of under-crossings
on $K_i$. This is the set of crossings such that the under-arcs belong to 
the component $K_i$. 
Then the state-sum
 $\Phi_i(L) =\sum_{{\cal C}}  \prod_{\tau \in UK_i}  B( \tau, {\cal C})$
 is an invariant of a link $L$ for each $i=1, \cdots, n$.
\end{lemma}

The vector 
$\vec{\Phi}(L)=( \Phi_i(L) )_{i=1}^n $ of the state-sum invariants is called
the {\it component-wise} (quandle) cocycle invariant of $L$.
We observe here that 
the theorem in this section applies to  component-wise cocycle invariants.

\begin{theorem} \label{linkcolorthm}
Let $\vec{\Phi}(L)=( \Phi_i(L) )_{i=1}^n $ be the component-wise 
cocycle invariant of a link $L=K_1 \cup \ldots \cup K_n$
with a quandle $X$ and a cocycle $\phi \in Z^2_{\rm Q}(X;A)$ for 
an abelian group $A$. 
Then  $\Phi_i(L)$ is not a positive integer for some $i$ if and 
only if there is a coloring of $L$ by $X$ that does not 
extend to a coloring of $L$ by $E(X, A, \phi)$. 
\end{theorem}

\begin{figure}
\begin{center}
\mbox{
\epsfxsize=3in
\epsfbox{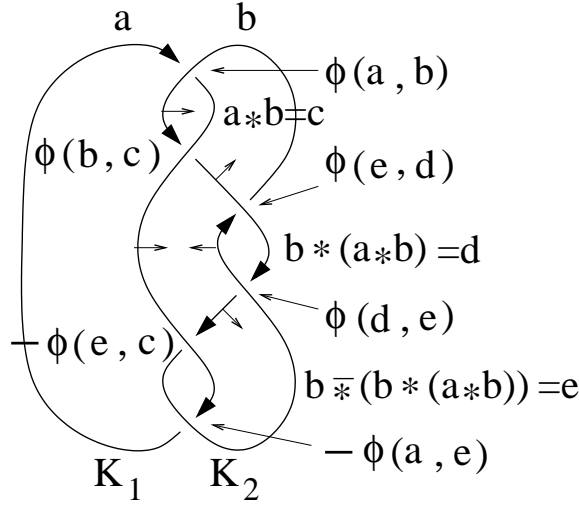} 
}
\end{center}
\caption{ A colored Whitehead link }
\label{whitehead} 
\end{figure}

\begin{example} {\rm 
In Fig.~\ref{whitehead}, 
a Whitehead link $L=K_1 \cup K_2$ is depicted. 
Let $\phi \in Z^2_{\rm Q}(R_8; \Z_2)$ be 
the
cocycle defined in
Cor.~\ref{dihedlemma}, that is, $\phi$ defines
the extension $E=R_{16}=E(R_8, \Z_2, \phi)$.
We evaluate the component-wise cocycle invariant
$\vec{\Phi}(L)=( \Phi_1(L),  \Phi_2(L) ) $.
Denote the generator of the coefficient group $\Z_2$ by $t$,
so that $\Z_2=\{ 1, t \}$ and the invariant takes the form 
of $\vec{\Phi}(L)=( A_1 + B_1 t ,  A_2 + B_2 t ) $,
where $A_i, B_i$ ($i=1,2$) are positive integers.

The colors assigned to arcs are represented by 
the letters
$a$ through $e$.
{}From the figure, 
it is seen that all the colors are determined by the colors $a$ and $b$
assigned to the top two arcs. 
It is seen by calculations that for any choice of two elements of $R_8$ 
for $a$ and $b$, there is a unique coloring of $L$ by $R_8$ that 
restricts to the chosen elements for $a$ and $b$. 
Therefore, there are $8^2=64$ colorings of $L$ by $R_8$. 

We show that the state-sum term 
$ \prod_{\tau \in UK_1}  B( \tau, {\cal C})$ is $1$ if and only if 
$a$ and $b$ have the same parity (both even or both odd). 

Suppose that $a$ and $b$ are 
both even, so that $a=2 \alpha$, $b=2 \beta$. 
Then one computes that $c=4 \beta - 2 \alpha$, $d=6 \beta - 4 \alpha$,
and we obtain $e=2 \beta=b$. 
Similar computations show that $e=b$ if $a$ and $b$ are both odd. 
{}From the figure, the state-sum term for $UK_1$ is 
$\phi(a, b) \phi(a, e)^{-1}$, which is equal to 
 $\phi(a, b) \phi(a, b)^{-1}=1$, in this case.
Suppose now that $a$ and $b$ have opposite parities. 
By setting $a=2 \alpha +1$ and $b = 2 \beta$ (and vice versa), 
we compute that $e=b+4$, so that we obtain 
the state-sum term $\phi(a, b) \phi(a, e)^{-1}=\phi(a, b) \phi(a, b+4)^{-1}$.
We claim that this is $t$. 

Using the formula at the end of the proof of Theorem~\ref{firstextthm},
we have $\phi(a,b)=[s(a)*s(b) - s(a*b)]/8$ (additively). 
Here, $s(a)*s(b)$ is $2b-a$ computed modulo $16$, and 
$s(a*b)$ is  $2b-a$ computed modulo $8$, then regarded as an element
modulo $16$. Since $a*(b+4)=2(b+4)-a=(2b-a)+8$ modulo $16$, we have
$$ \phi(a,b)-\phi(a, b+4) = [s(a)*s(b) - s(a*b)]/8 -   [s(a)*s(b)+8 - s(a*b)]/8
= 1 \; \mbox{mod}\; (2)  $$
written additively. This proves the above claim. 
There are $32$ colorings with the same parity, and $32$ with distinct parities.
Hence we obtain 
 $\vec{\Phi}(L)=( 32 + 32t, 32+ 32t ) $.

Theorem~\ref{linkcolorthm} implies that there are colorings by $R_8$ that 
do not extend to colorings by $R_{16}$. 
In fact, from the proof of Theorem~\ref{knotcolorthm}, we see that 
$32$ colorings having the same parity for $a$ and $b$ extend to
$R_{16}$, and those $32$ colorings  with the opposite parities do not.
This fact can be computed directly, and gives an alternate method of computing
the above invariant using Cor.~\ref{z2case}. 

} \end{example}

\noindent
{\large\bf Acknowledgments} \
JSC is being supported by NSF grant DMS-9988107.
MS is being supported by NSF grant DMS-9988101.

\end{document}